\documentclass[12pt,a4paper,reqno]{article}    
\usepackage{graphics}
\usepackage{cite}
\usepackage[demo]{graphicx}
\usepackage[fleqn]{amsmath}
\usepackage{booktabs}
\usepackage{multirow}
\usepackage{epsfig}
\usepackage{epstopdf}
\usepackage{booktabs}
\usepackage{placeins}
\usepackage{anysize}
\usepackage{enumerate}
\usepackage{latexsym,amssymb,amsmath,amscd,amsthm,amsxtra}
\usepackage{subfigure}
\usepackage{tabularx}
\usepackage[table]{xcolor}
\allowdisplaybreaks
\definecolor{lightgray}{gray}{0.9}
\newtheorem{theorem}{Theorem}
\newtheorem{prop}{Proposition}

\numberwithin{equation}{section}

\def\({\left( }
\def\){\right )}
\begin{document}
\title {\textbf{The application of cubic trigonometric B-spline to the numerical solution of time-fractional telegraph equation}}
\author{Muhammad Yaseen, Muhammad Abbas\thanks{Corresponding authors: m.abbas@uos.edu.pk}\\
{Department of Mathematics, University of Sargodha, Pakistan.}\\
\\
}
\date{}
\maketitle
\begin{abstract}
In this paper,  an efficient numerical technique for the time-fractional telegraph equation is proposed.  The aim of this paper is to use a relatively new type of B-spline called the  cubic trigonometric B-splines for the proposed scheme.  This technique is based on finite difference formulation for the Caputo time-fractional derivative and cubic trigonometric B-splines based technique for the derivatives in space. A stability analysis of the scheme is set up to affirm that the errors do not amplify.  Computational  experiments are carried out in addition to  verify the theoretical analysis.  Numerical results are compared with some existing techniques and it is concluded that the present scheme is more accurate and effective.  
\end{abstract}
\maketitle
\section{Introduction}
In recent years, the tools of fractional calculus have been successfully used to describe many physical phenomena in  science and engineering \cite{pod,hil,kilbas}. Recently, there have been reporting of many applications typically expressed by fractional partial differential equations (FPDEs). The importance of FPDEs lies in the fact that the solutions offered by FPDEs have descriptions that well approximate the chemical, physical and biological phenomena than their integer order counterparts. As a result, FPDEs have attained special status among scientists and engineers.

A number of phenomenon such as propagation of electric signals \cite{jordan}, transport of neutron in a nuclear reactor \cite{nata} and random walks \cite{mika} are described by a class of hyperbolic partial differential equations called the fractional telegraph equations \cite{hatch}. The general form of the time-fractional telegraph equation is given by
\begin{equation}\label{1.1}
\begin{cases}
    \frac{\partial^ {\gamma} u(x,t)}{\partial t^{\gamma}}+\gamma_1 \frac{\partial^ {\gamma-1} u(x,t)}{\partial t^{\gamma-1}} +\gamma_2 u(x,t)=\gamma_3 \frac{\partial^2 u(x,t)}{\partial x^2}+f(x,t),\\
    \qquad\qquad\qquad\qquad a\leq x \leq b,\,\,\,1 < \gamma < 2,\,\,\,0\leq t\leq T.
    \end{cases}
\end{equation}
with initial conditions
\begin{equation}\label{1.2}
\left\{
\begin{array}{cc}
 u(x,0)=\phi_1(x),\\
 u_t(x,0)=\phi_2(x)
\end{array}\right.
\end{equation}
and the boundary conditions
\begin{equation}\label{1.3}
\left\{
\begin{array}{cc}
 u(a,t)=\psi_1(t),\\
 u(b,t)=\psi_2(t)
\end{array}\qquad\qquad\qquad t\geq 0,\right.
\end{equation}
where $a,b, \phi_1(x),\phi_2(x), \psi_1(t)$ and $\psi_2(t)$ are given and  $\frac{\partial^\gamma}{\partial t^\gamma}u(x,t)$ represents the Caputo fractional derivative of order $\gamma$  given by \cite{pod,hil,kilbas}
\begin{equation}\label{1.4}
\frac{\partial^\gamma}{\partial t^\gamma}u(x,t)=
\begin{cases}
\frac{1}{\Gamma(n-\gamma)}\int\limits_0^t \frac{\partial^n u(x,s)}{\partial s^n}(t-s)^{n-\gamma-1}ds,\qquad n-1<\gamma <n\\
\frac{\partial^n u(x,s)}{\partial s^n},\qquad\qquad\qquad\qquad\qquad\qquad \gamma=n.
\end{cases}
\end{equation}
In addition, $\gamma_1, \gamma_2, \gamma_3$ are given constants. Note that in case of $\gamma=2$, Eq. (\ref{1.1}) corresponds to the classical second-order telegraph equation.

Various numerical and analytical methods are accessible in literature for the time-fractional telegraph equation. Tasbozan and Esen \cite{esen1} utilized B-spline Galerkin method for the numerical solutions fractional telegraph equation. Hosseini \textit{et al.}\cite{vrh} have made use of radial basis functions to obtain the numerical solution of time-fractional telegraph equation. Sweilam \textit{et al.}\cite{swlam} have used Sinc-Legendre collocation procedure to find an approximate solution of time-fractional-order telegraph equation. A classic work of  Orsinger and Zhao \cite{sing} regarding the space-fractional telegraph equation and the related fractional telegraph process appeared in 2003.   S. Momani \cite{mom1} obtained analytic and approximate solutions of the space and time-fractional telegraph equations. Chen \textit{et al.} \cite{chn1} utilized the method of separating variables to obtain analytical solutions for the time-fractional telegraph equation. Wei \textit{et al.}\cite{wei} presented a fully discrete local discontinuous Galerkin method for solving the fractional telegraph equation. Wang \textit{et al.} \cite{wan1} used  reproducing kernel for solving a class of time-fractional telegraph equation with initial value conditions. Hashemi and Baleanu \cite{hash} ]utilized a geometric approach and the method of lines to obtain a numerical approximation of higher-order time-fractional telegraph equation. Jiang and Lin \cite{lin1} obtained  the exact solution of the time-fractional telegraph equation in the reproducing kernel space.  In \cite{kumar}, Kumar presented a  new analytical modeling for fractional-telegraph equation via Laplace transform. Mollahasani \textit{et al.} \cite{mola} developed a new technique based on hybrid functions for the numerical treatment of telegraph equations of fractional order. Hariharan  \textit{et al.} \cite{hari} utilized a wavelet method for a class of space and time-fractional telegraph equations. Analytical solutions of space and time-fractional telegraph equations were obtained by Yildirim \cite{yal} by using He's homotopy perturbation method.

The principal purpose of this paper is to present a numerical scheme for the time-fractional telegraph equation that is computationally efficient and provides better results than some existing  numerical procedures \cite{esen1,vrh,swlam}. To authors knowledge this paper is first attempt towards finding the numerical solution of time-fractional telegraph equations using cubic trigonometric  B-splines. A detailed stability analysis of the scheme is presented to assert that errors do not amplify. Numerical experiments are performed to further set the accuracy and validity of the technique.

The rest of the paper is organized as follows. In section 2, the numerical scheme primarily based on cubic trigonometric B-splines is derived in detail. Section 3 discusses the stability analysis. Section 4 indicates a comparison of our numerical consequences with those of \cite{esen1,vrh,swlam}. Section 5 summarizes the conclusions of this study.
\section{The Derivation of the Scheme}
For given positive integers $M$ and $N$, let $\tau=\frac{T}{N}$ be the temporal and $h=\frac{b-a}{M}$ the spatial step sizes respectively. Following the usual notations, set $t_n=n\tau~(0\leq n \leq N)$, $x_j=jh,~(0 \leq j \leq M)$, $\Omega_{\tau}=\{t_n|0 \leq n \leq N\}$ and $\Omega_h=\{x_j|0 \leq j \leq M\}$. Let $u_j^n$ be approximation to exact solution at the point $(x_j,t_n)$ and $\mathcal{A}=\{u_j^n|0\leq j\leq M, 0\leq n\leq N\}$ be grid function space defined on $\Omega_h \times \Omega_\tau$.  The solution domain $a\leq x\leq b$ is uniformly partitioned  by knots $x_{i}$ into $N$ subintervals $[x_{i}, x_{i+1}]$ of equal length $h$, $i=0,1,2,...,M-1$, where $a=x_{0}<x_{1}<...<x_{n-1}<x_{M}=b$. Our scheme for solving (\ref{1.1}) requires approximate solution $U(x,t)$ to the exact solution $u(x,t)$  in the following form  \cite{36,37}
\begin{equation}\label{2.1}
    U(x,t)=\sum\limits_{i=-1}^{N-1} c_i(t) TB_i^4 (x),
\end{equation}
where $c_i(t)$ are unknowns to be determined and $TB^4_i (x)$  \cite{abbas} are twice differentiable cubic Trigonometric basis functions given by
\begin{equation}\label{2.2}
    TB^4_i (x)=\frac{1}{w}
    \begin{cases}
    p^3(x_i) & x\in [x_i,x_{i+1}]\\
    p(x_i)(p(x_i)q(x_{i+2})+q(x_{i+3})p(x_{i+1}))+q(x_{i+4})p^2(x_{i+1}),  & x \in [x_{i+1},x_{i+2}]\\
    q(x_{i+4})(p(x_{i+1})q(x_{i+3})+q(x_{i+4})p(x_{i+2}))+p(x_i)q^2(x_{i+3}), & x \in [x_{i+2},x_{i+3}]\\
    q^3(x_{i+4}), & x \in [x_{i+3},x_{i+4}]
    \end{cases}
\end{equation}
where\\
$\displaystyle{p\(x_i\)=\sin\(\frac{x-x_i}{2}\), q\(x_i\)=\sin\(\frac{x_i-x}{2}\), w=\sin\(\frac{h}{2}\) \sin\(h\) \sin\(\frac{3 h}{2}\).}$\\
 Due to local support property of the cubic trigonometric B-splines only $TB^4_{j-1}(x), TB^4_{j}(x)$ and $TB^4_{j+1}(x)$ are survived so that the approximation $u_j^n$ at the grid point $(x_j,t_n)$ at $n^{th}$  time level is given as:
\begin{equation}\label{2.3}
    u(x_{j},t_{n})=u_j^n=\sum\limits_{j=i-1}^{i+1} c_j^n(t) TB^4_j (x).
\end{equation}
The  time dependent unknowns $c_j^n(t)$ are to be determined by making use of the initial and boundary conditions, and the collocation conditions on $TB^4_i (x)$. As a result the approximations $u_j^n$ and its necessary derivatives are given as:
\begin{equation}\label{2.4}
    \begin{cases}
   \displaystyle{ u_j^n=a_1 c_{j-1}^n+a_2 c_{j}^n+ a_1 c_{j+1}^n}\\
    \displaystyle{ (u_j^n)_x=-a_3 c_{j-1}^n+a_3 c_{j+1}^n}\\
    \displaystyle{ (u_j^n)_{xx}=a_4 c_{j-1}^n+ a_5 c_{j}^n +a_4 c_{j+1}^n,}
    \end{cases}
\end{equation}
where\\
$\displaystyle{a_1=\csc\(h\) \csc\(\frac{3h}{2}\)\sin^2\(\frac{h}{2}\),}$\\
$\displaystyle{a_2=\frac{2}{1+2 \cos\(h\)},}$\\
$\displaystyle{a_3=\frac{3}{4} \csc \(\frac{3h}{2}\),}$\\
$\displaystyle{a_4=\frac{3+9\cos\(h\)}{4 \cos\(\frac{h}{2}\)-4 \cos\(\frac{5 h}{2}\)},}$\\
$\displaystyle{a_5=-\frac{3 \cot^2\(\frac{h}{2}\)}{2+4 \cos\(h\)}.}$\\\\
Following \cite{vrh}, the fractional derivatives  $\frac{\partial^\gamma}{\partial t^\gamma}u(x,t)$  and $\frac{\partial^{\gamma-1}}{\partial t^{\gamma-1}}u(x,t)$ are  discrtetized as:
\begin{equation}\label{2.5}
    \frac{\partial^\gamma}{\partial t^\gamma}u(x_j,t_{n+1})=\alpha_0 \sum\limits_{k=0}^{n} b_k (u^{n-k+1}-2u^{n-k}+u^{n-k-1}),\qquad 1<\gamma <2
\end{equation}
and
\begin{equation}\label{2.6}
    \frac{\partial^{\gamma-1}}{\partial t^{\gamma-1}}u(x_j,t_{n+1})=\alpha_0\Delta t \sum\limits_{k=0}^{n} b_k (u^{n-k+1}-u^{n-k}), \qquad 0< \gamma-1 < 1,
\end{equation}
where $\alpha_0=\frac{\Delta t^{-\gamma}}{\Gamma[3-\gamma]}$ and  $b_k=(k+1)^{1-\gamma}-k^{1-\gamma}$. It is straight forward to confirm that
\begin{itemize}
  \item $b_k>0,~~k=0,1,\cdots ,n$.
  \item $1=b_0>b_1>b_2>\cdots> b_n~\text{and}~b_n\rightarrow 0~\text{as}~n\rightarrow \infty$.
  \item $\sum\limits_{k=0}^n (b_k-b_{k+1})=1$.
\end{itemize}
To obtain temporal discretization, we substitute (\ref{2.5}) and (\ref{2.6}) into (\ref{1.1}) to get:
\begin{equation}\label{2.7}
\begin{split}
\alpha_0(u^{n+1}-2u^{n}+u^{n-1})+\alpha_0\sum\limits_{k=1}^n b_k(u^{n+1-k}-2u^{n-k}+u^{n-1-k})+\gamma_1\Delta t\alpha_0 (u^{n+1}-u^{n})\\+\gamma_1 \Delta t \alpha_0 \sum\limits_{k=1}^n b_k(u^{n+1-k}-u^{n-k})+\gamma_2 u^{n+1}-\gamma_3 \frac{\partial^2 u^{n+1}}{\partial x^2}=f_j^{n+1}.
\end{split}
\end{equation}
It is observed that the term $u^{-1}$ will appear when $n=0$ or $k=n$. Using the central forward difference formula, we utilize the given initial condition to obtain
\begin{equation}\label{2.8}
u_t^0=\frac{u^1-u^{-1}}{2\Delta t}
\end{equation}
from where we observe that $u^{-1}=u^1-2\Delta t \phi_2(x).$\\
To obtain full discretization, we substitute the approximations (\ref{2.4}) into (\ref{2.7}) and get
\begin{eqnarray}\label{2.9}
&&((\alpha_0+\gamma_1\Delta t \alpha_0+\gamma_2)a_1-\gamma_3 a_4)c_{j-1}^{n+1}+((\alpha_0+\gamma_1\Delta t \alpha_0+\gamma_2)a_2-\gamma_3 a_5)c_{j}^{n+1}\nonumber\\&&+((\alpha_0+\gamma_1\Delta t \alpha_0+\gamma_2)a_1-\gamma_3 a_4)c_{j}^{n+1}\nonumber\\&=& (2\alpha_0+\gamma_1\Delta t \alpha_0)(a_1 c_{j-1}^n+a_2 c_{j}^n+a_1 c_{j+1}^n)-\alpha_0 (a_1 c_{j-1}^{n-1}+a_2 c_{j}^{n-1}+a_1 c_{j+1}^{n-1})\nonumber\\&&- \gamma_1\alpha_0\Delta t\sum\limits_{k=1}^n b_k\left(a_1(c_{j-1}^{n+1-k}-c_{j-1}^{n-k})+a_2(c_{j}^{n+1-k}-c_{j}^{n-k})+a_1(c_{j+1}^{n+1-k}-c_{j+1}^{n-k})\right)\nonumber\\&&-\alpha_0 \sum\limits_{k=1}^n b_k(a_1(c_{j-1}^{n+1-k}-2 c_{j-1}^{n-k}+c_{j-1}^{n-1-k})+a_2(c_{j}^{n+1-k}-2 c_{j}^{n-k}+c_{j}^{n-1-k})\nonumber\\&&+a_1(c_{j+1}^{n+1-k}-2 c_{j+1}^{n-k}+c_{j+1}^{n-1-k}))+f_j^{n+1}.
\end{eqnarray}
The equation (\ref{2.9}) consists of $(N-1)$ linear equations in $N+1$ unknowns. To obtain a unique solution to the system, we need two additional equations which can be obtained by utilizing the given boundary conditions (\ref{1.3}). As a result a diagonal matrix of dimension $(N+1)\times (N+1)$ is obtained which can be solved by using any suitable numerical algorithm.
\section{Stability Analysis}
This section deals with the stability analysis of the fully discrete scheme (\ref{2.9}). By Duhamels' principle \cite{dhuml} it can be concluded that the stability analysis for an inhomogeneous problem is a direct outcome of the analysis for the corresponding homogeneous case. So it is sufficient to present the stability analysis for the force free case $f=0$. In this study, we assume the growth factor of a Fourier mode to be $\sigma_j^k$ and let $\tilde{\sigma}_j^k$ be its approximation. Define $E_j^k=\sigma_j^k-\tilde{\sigma}_j^k$ so that from (\ref{2.9}), we obtain the following round off error equation
\begin{eqnarray}\label{3.1}
&&((\alpha_0+\gamma_1\Delta t \alpha_0+\gamma_2)a_1-\gamma_3 a_4)E_{j-1}^{n+1}+((\alpha_0+\gamma_1\Delta t \alpha_0+\gamma_2)a_2-\gamma_3 a_5)E_{j}^{n+1}\nonumber\\&&+((\alpha_0+\gamma_1\Delta t \alpha_0+\gamma_2)a_1-\gamma_3 a_4)E_{j}^{n+1}\nonumber\\&=& (2\alpha_0+\gamma_1\Delta t \alpha_0)(a_1 E_{j-1}^n+a_2 E_{j}^n+a_1 E_{j+1}^n)-\alpha_0 (a_1 E_{j-1}^{n-1}+a_2 E_{j}^{n-1}+a_1 E_{j+1}^{n-1})\nonumber\\ &&-\gamma_1\alpha_0\Delta t\sum\limits_{k=1}^n b_k\left(a_1(E_{j-1}^{n+1-k}-E_{j-1}^{n-k})+a_2(E_{j}^{n+1-k}-E_{j}^{n-k})+a_1(E_{j+1}^{n+1-k}-E_{j+1}^{n-k})\right)\nonumber\\&&-\alpha_0 \sum\limits_{k=1}^n b_k(a_1(E_{j-1}^{n+1-k}-2 E_{j-1}^{n-k}+E_{j-1}^{n-1-k})+a_2(E_{j}^{n+1-k}-2 E_{j}^{n-k}+E_{j}^{n-1-k})\nonumber\\&&+a_1(E_{j+1}^{n+1-k}-2 E_{j+1}^{n-k}+E_{j+1}^{n-1-k})).
\end{eqnarray}
The error equation satisfies the boundary conditions
\begin{equation}\label{3.2}
    E_0^k=\psi_1(t_k),~~E_M^k=\psi_2(t_k),~~k=0,1,\cdots, N
\end{equation}
and the initial conditions
\begin{equation}\label{3.3}
    E_j^0=\phi_1(x_j),~~(E_t)_j^0=\phi_2(x_j),~~j=1,2,\cdots, M.
\end{equation}
Define the grid function
$$E^k(x)=\begin{cases}
      E_j^k, & \qquad x_j-\frac{h}{2} < x \leq x_j+\frac{h}{2} j=1,\cdots,M-1 \\
      0, &  \qquad a<x\leq \frac{h}{2} \text~{or}~ b-\frac{h}{2}<x\leq b.
   \end{cases}$$
Note that the Fourier expansion of $E^K(x)$ is
$$E^k(x)=\sum\limits_{m=-\infty}^\infty a_k(m)e^{\frac{i 2 \pi m x}{(b-a)}}, k=0,1,\cdots,N$$
where $a^k(m)=\frac{1}{(b-a)}\int\limits_a^b E^k(x)e^{\frac{-i 2 \pi m x}{(b-a)}} dx$. Let $$E^k=[E_1^k,E_2^k,\cdots,E_{M-1}^k]^T$$ and introduce the norm: $$\|E^k\|_2=\left(\sum\limits_{j=1}^{M-1} h |E_j^k|^2\right)^{\frac{1}{2}}=\left[\int\limits_a^b |E^k(x)|^2 dx\right]^{\frac{1}{2}}.$$
By  Parseval equality, it is observed that
$$\int\limits_a^b |E^k(x)|^2 dx=\sum\limits_{m=-\infty}^{\infty} |a_k(m)|^2,$$ so that the following relation is obtained
\begin{equation}\label{3.4}
\|E^k\|_2^2=\sum\limits_{m=-\infty}^\infty |a_k(m)|^2.
\end{equation}
Suppose that equations (\ref{3.1})-(\ref{3.3}) have solution of the form $E_j^n=\xi_n e^{i\beta j h}$, where $i=\sqrt{-1}$ and $\beta$ is real. Substituting this expression into (\ref{3.1}), dividing by $e^{i\beta j h}$, using the relation $e^{-i\beta h}+e^{i\beta h}=2cos(\beta h)$ and collecting the like terms, we obtain
\begin{eqnarray}\label{3.5}
&&\left(((\alpha_0+\gamma_1\Delta t \alpha_0+\gamma_2)a_1-\gamma_3 a_4)2\cos(\beta h)+((\alpha_0+\gamma_1\Delta t \alpha_0+\gamma_2)a_2-\gamma_3 a_5)\right)\xi_{n+1}\nonumber\\&&\nonumber \\&=& \left((2\alpha_0+\gamma_1\Delta t \alpha_0)a_1 2\cos(\beta h)+(2\alpha_0+\gamma_1\Delta t \alpha_0)a_2 \right)\xi_n-(\alpha_0 a_1 2\cos(\beta h)+\alpha_0 a_2 )\xi_{n-1}\nonumber\\ &&-\gamma_1\alpha_0\Delta t\sum\limits_{k=1}^n b_k\left((2 a_1\cos(\beta h)+a_2)\xi_{n+1-k}-(2 a_1\cos(\beta h)+a_2)\xi_{n-k}\right)\nonumber\\&&-\alpha_0 \sum\limits_{k=1}^n b_k((2a_1 \cos(\beta h)+a_2)\xi_{n+1-k}-2(2a_1 \cos(\beta h)+a_2)\xi_{n-k}\nonumber\\&&+(2a_1 \cos(\beta h)+a_2)\xi_{n-1-k})).\nonumber\\&&
\end{eqnarray}
Without loss of generality, we can assume that $\beta=0$, so that (\ref{3.5})reduces to
\begin{equation}\label{3.6}
\begin{split}
\xi_{n+1}=\frac{2+\gamma_1 \Delta t}{\nu}\xi_n-\frac{1}{\nu}\xi_{n-1}-\frac{\gamma_1 \Delta t}{\nu}\sum\limits_{k=1}^n b_k (\xi_{n+1-k}-\xi_{n-k})\\-\frac{1}{\nu}\sum\limits_{k=1}^n b_k (\xi_{n+1-k}-2\xi_{n-k}+\xi_{n-1-k})
\end{split}
\end{equation}
where $\nu=\left( 1+\gamma_1 \Delta t+\frac{\gamma_2}{\alpha_0}-\frac{\gamma_3}{\alpha_0}(\frac{2a_4+a_5}{2a_1+a_2})\right).$
\begin{prop}\label{aa}
If $\xi_n$ is the solution of equation (\ref{3.6}), then $|\xi_n|\leq 2|\xi_0|$, $n=0,1,\cdots, T \times N$ provided $v\geq 2+\gamma_1 \Delta t$.
\end{prop}
\textbf{Proof}. Mathematical induction is used to prove the result. For $n=0$, we have from equation $(\ref{3.6})$ that $\xi_1=\frac{2+\gamma_1\Delta t}{\nu} \xi_0$ and since $|\frac{2+\gamma_1\Delta t}{\nu}|\leq 1$, therefore
$$|\xi_1|= |\frac{2+\gamma_1\Delta t}{\nu}||\xi_0|\leq |\xi_0|\leq 2 |\xi_0|. $$
Now suppose that $|\xi_n|\leq 2|\xi_0|$, $n=1,\cdots,T \times N-1 $ so that from (\ref{3.6}), we obtain
\begin{eqnarray}
  |\xi_{n+1}|&\leq& |\xi_n|+\frac{|\xi_{n-1}|}{|\nu|}+|\frac{\gamma_1 \Delta t}{\nu}|\sum\limits_{k=1}^n b_k (|\xi_{n+1-k}|-|\xi_{n-k}|)\nonumber\\&&+|\frac{1}{\nu}|\sum\limits_{k=1}^n b_k (|\xi_{n+1-k}|-2|\xi_{n-k}|+|\xi_{n-1-k}|)\\ \nonumber
    &\leq& |\xi_0|+|\xi_0|+|\frac{\gamma_1 \Delta t}{\nu}|\sum\limits_{k=1}^n b_k (|\xi_{0}|-|\xi_{0}|)\nonumber\\&&+|\frac{1}{\nu}|\sum\limits_{k=1}^n b_k (|\xi_{0}|-2|\xi_{0}|+|\xi_{0}|)\\ \nonumber
        &=& 2|\xi_0|.
\end{eqnarray}
This completes the proof.
\begin{theorem}
The collocation scheme (\ref{3.5}) is stable provided that $v\geq 2+\gamma_1 \Delta t$.
\end{theorem}
\noindent
\textbf{Proof}. Using formula (\ref{3.4}) and Proposition \ref{aa}, we  obtain $$\|E^k\|_2\leq 2\|E^0\|_2, k=0,1,\cdots N $$  which establishes that the  scheme is stable under the given condition.

\section{Concluding Remarks}\label{s6}
This study presents a  numerical technique based on cubic trigonometric B-spline for the time-fractional telegraph equation. The scheme utilized the usual finite difference scheme to approximate the Caputo time-fractional derivative and the derivative in space are approximated using the cubic trigonometric B-spline basis functions. A special attention has been given to study the stability analysis of the scheme. The obtained results are compared with those of some existing techniques. The comparison reveals that the presented scheme is comparable with other existing techniques for time-fractional telegraph equation   in terms of accuracy, flexibility and efficiency. Moreover, the scheme can be applied to a large class of  fractional order partial differential equations.

\end{document}